\documentclass[11pt]{article}
\usepackage{amsfonts}
\usepackage{latexsym, amssymb, amsmath, amscd, amsfonts, epsfig, graphicx}
\newtheorem{thm}{Theorem}[section]

\def\pf{\noindent{\it Proof.} }
\def\qed{\nopagebreak\hfill{\rule{4pt}{7pt}}
\medbreak}

\def\qed{\nopagebreak\hfill{\rule{4pt}{7pt}}
\medbreak}

\makeatletter
\def\ExtendSymbol#1#2#3#4#5{\ext@arrow 0099{\arrowfill@#1#2#3}{#4}{#5}}
\makeatother

\title{An Overpartition Analogue of Bressoud's Theorem of Rogers-Ramanujan Type}
\author{William Y.C. Chen\raisebox{5pt}{\scriptsize 1},
Doris D. M. Sang\raisebox{5pt}{\scriptsize 2}, and Diane Y. H.
Shi\raisebox{5pt}{\scriptsize 3}}
\date{Center for Combinatorics, LPMC-TJKLC\\
 Nankai University\\
Tianjin 300071, P.R. China \\
\vspace{15pt} \raisebox{5pt}{\scriptsize 1\,}chen@nankai.edu.cn,
\raisebox{5pt}{\scriptsize 2\,}sdm@cfc.nankai.edu.cn,
\raisebox{5pt}{\scriptsize 3\,}yahuishi@gmail.com}

\begin{document}
\maketitle
\noindent {\bf Abstract.} For $k\geq i\geq 1$, let $B_{k,i}(n)$ denote the number of partitions of $n$ such that part $1$ appears at most $i-1$ times,  two consecutive integers
$l$ and $l+1$ appear at most $k-1$ times and if  $l$ and $l+1$ appear exactly $k-1$ times then the total sum of the parts $l$ and $l+1$ is congruent to $i-1$ modulo $2$. Let $A_{k,i}(n)$ denote the number of partitions with  parts not congruent to  $i$, $2k-i$ and $2k$ modulo $2k$. Bressoud's theorem states that $A_{k,i}(n)=B_{k,i}(n)$.  Corteel,  Lovejoy, and  Mallet
found an overpartition analogue of Bressoud's theorem for   $i=1$, that is, for partitions not containing nonoverlined part $1$.
We obtain
  an overpartition analogue of Bressoud's theorem in the general case.
For $k\geq i\geq 1$, let $D_{k,i}(n)$ denote the number of overpartitions of $n$ such that the nonoverlined part $1$ appears at most $i-1$ times,  for any integer $l$, $l$ and nonoverlined $l+1$ appear at most $k-1$ times and if  the parts $l$ and the nonoverlined part $l+1$ appear exactly $k-1$ times then the total sum of the parts $l$ and nonoverlined part $l+1$ is congruent to the number of overlined parts that are less than $l+1$ plus $i-1$ modulo $2$. Let $C_{k,i}(n)$ denote the number of overpartitions with the nonoverlined parts not congruent to  $\pm i$ and $2k-1$ modulo $2k-1$.
 We show that $C_{k,i}(n)=D_{k,i}(n)$. This relation can also be considered as
 a Rogers-Ramanujan-Gordon  type theorem for overpartitions.

\noindent {\bf Keywords:} Rogers-Ramanujan-Gordon theorem, overpartition, Bressoud's theorem,

\noindent {\bf AMS Subject Classification:} 05A17, 11P84

\section{ Introduction}

The Rogers-Ramanujan-Gordon theorem is a combinatorial generalization of the Rogers-Ramanujan identities \cite{ram19,rog1894}, see Gordon \cite{gor61}. It establishes the equality between
the number of partitions of $n$ with parts satisfying  certain residue
 conditions and the number of partitions of $n$ with certain difference conditions. Gordon found an involution for an equivalent form of the
  generating function identity for this relation.
  An algebraic proof was given by Andrews \cite{and66} by using a recursive approach.
   It should be noted that the Rogers-Ramanujan-Gordon theorem is concerned only
   with odd moduli. Bressoud \cite{Bre79} succeeded in finding a theorem of Rogers-Ramanujan-Gordon type for even moduli by using an  algebraic approach
    in the spirit of Andrews \cite{and76}.

The objective of this paper is to give an overpartition analogue of
Bressoud's theorem. We shall derive  the equality between the number of overpartitions of $n$ such that the nonoverlined parts belong to certain residue classes modulo some odd positive integer  and the number of overpartitions of $n$ with parts satisfying certain  difference conditions.   A special case of this relation has been discovered by Corteel,  Lovejoy, and  Mallet \cite{CLM}.

An  overpartition analogue  of the Rogers-Ramanujan-Gordon theorem was obtained by Chen, Sang and Shi \cite{chen},  which states that the number of overpartitions of $n$ with  nonoverlined parts  belonging to certain residue classes
modulo some even  positive integer  equals  the number of overpartitions  of $n$ with parts satisfying certain difference conditions. However, as will be seen,
 the proof of the overpartition analogue of the Rogers-Ramanujan-Gordon theorem does not seem to be directly applicable to the
case for the overpartition analogue of Bressoud's theorem.

Let us give an overview of  some definitions.
A partition $\lambda$ of a positive integer $n$ is a non-increasing sequence of positive integers $\lambda_1\geq \cdots\geq \lambda_s>0$ such that $n=\lambda_1+\cdots+\lambda_s$. The partition of zero is the partition with no parts.
An overpartition $\lambda$ of a positive integer $n$ is also a non-increasing sequence of positive integers $\lambda_1\geq \cdots\geq \lambda_s>0$ such that $n=\lambda_1+\cdots+\lambda_s$ and the first occurrence  of  each integer may be overlined. For example, $(\overline{7},7,6,\overline{5},2,\overline{1})$ is an overpartition of $28$.  Many $q$-series identities  have combinatorial
interpretations in terms of overpartitions,
see, for example, Corteel and Lovejoy \cite{cor04}. Furthermore, overpartitions possess many analogous properties of ordinary partitions, see Lovejoy \cite{lov03, lov06}. For example,
various   overpartition theorems  of the Rogers-Ramanujan-Gordon type have been
obtained by  Corteel and
Lovejoy \cite{cor07}, Corteel, Lovejoy and Mallet \cite{CLM} and Lovejoy \cite{lov03, lov04, lov07}.
For a partition or an overpartition $\lambda$ and for any integer $l$, let $f_l(\lambda) (f_{\overline{l}}(\lambda))$ denote the number of occurrences of $l$
non-overlined (overlined) in $\lambda$. Let $V_{\lambda}(l)$ denote the number of overlined parts in $\lambda$ that are less than or
equal to $l$.

We shall adopt  the common notation  as used in Andrews \cite{and76}. Let
 \[(a)_\infty=(a;q)_\infty=\prod_{i=0}^{\infty}(1-aq^i),\]
 and \[(a)_n=(a;q)_n=\frac{(a)_\infty}{(aq^n)_\infty}.\]
We also write
 \[(a_1,\ldots,a_k;q)_\infty=(a_1;q)_\infty\cdots(a_k;q)_\infty.\]

The Rogers-Ramanujan-Gordon theorem reads as follows.

\begin{thm}(Rogers-Ramanujan-Gordon)\label{Gordon} For $k\geq i\geq 1$,
let $F_{k,i}(n)$ denote the number of partitions
of $n$ of the form $\lambda_1 + \lambda_2 + \cdots + \lambda_s$, where $\lambda_j \geq \lambda_{j+1}$,
$\lambda_j-\lambda_{j+k-1}\geq2$ and part $1$ appears at most $i-1$ times.
 Let $E_{k,i}(n)$ denote the number of
partitions of $n$ into parts $\not \equiv0,\pm i\; (mod\ 2k + 1)$.
Then for any $n \geq0$, we have
 \begin{equation}\label{ef}
 E_{k,i}(n) = F_{k,i}(n).
 \end{equation}
\end{thm}

In the algebraic proof of the above relation,  Andrews \cite{and66, and76} introduced  a hypergeometric function $J_{k,i}(a;x;q)$ as given by
\begin{equation} \label{jki}
J_{k,i}(a;x;q)=
H_{k,i}(a;xq;q)-axqH_{k,i-1}(a;xq;q),
\end{equation}
where
\begin{equation}\label{eqH}
H_{k,i}(a;x;q)=
\sum_{n=0}^{\infty}\frac{x^{kn}q^{kn^2+n-in}a^n(1-x^iq^{2ni})
(axq^{n+1})_{\infty}(1/a)_n}{(q)_n(xq^n)_\infty}.
\end{equation}

To prove (\ref{ef}), Andrews considered a refinement of $F_{k,i}(n)$, that is, the number of partitions enumerated by $F_{k,i}(n)$ with exactly $m$ parts,  denoted by $F_{k,i}(m,n)$, and he showed that $J_{k,i}(-1/q;x;q)$ and the generating function of $F_{k,i}(m,n)$ satisfy the same recurrence relation with the same initial values.  Setting $x=1$ and using   Jacobi's triple product identity,  we find that $J_{k,i}(-1/q;1;q;)$ equals the generating function for $E_{k,i}(n)$. This yields that $E_{k,i}(n)=F_{k,i}(n)$.

The following Rogers-Ramanujan-Gordon type theorem for even moduli is due to Bressoud \cite{Bre79}.

\begin{thm}\label{Bressoud}For $k\geq i\geq 1$, let $B_{k,i}(n)$ denote the number of partitions of $n$ of the form $\lambda=\lambda_1+\lambda_2+\cdots+\lambda_s$ such that
(i) $f_1(\lambda)\leq i-1$, (ii) $f_l(\lambda)+f_{l+1}(\lambda)\leq k-1$, and (iii) if $f_l(\lambda)+f_{l+1}(\lambda)=k-1$, then $lf_l(\lambda)+(l+1)f_{l+1}(\lambda)  \equiv i-1 \;(\rm{mod}\; 2)$. Let $A_{k,i}(n)$ denote the number of
partitions of $n$ with parts not congruent to
$0,\pm i$ modulo $2k$. Then we have
\begin{equation}
A_{k,i}(n)=B_{k,i}(n).
\end{equation}
\end{thm}

The proof of Bressoud  also uses the hypergeometric function $J_{k,i}(-1/q;x;q)$. But he needs a recurrence relation for  $(-xq)_\infty  J_{(k-1)/2,i/2}(a;x^2;q^2)$.

Lovejoy \cite{lov03} found the following overpartition analogues of  Rogers-Ramanujan-Gordon theorem for the cases $i=1$ and $i=k$.

\begin{thm}\label{lov1}
For $k\geq 1$, let $\overline{B}_{k}(n)$ denote the number of
overpartitions of $n$ of the form $\lambda_1+\lambda_2+\cdots+\lambda_s$ such that
$\lambda_j-\lambda_{j+k-1}\geq1$ if $\lambda_j$ is overlined and $\lambda_j-\lambda_{j+k-1}\geq2$
otherwise. Let $\overline{A}_{k}(n)$ denote the number of overpartitions of
$n$ into parts not divisible by $k$. Then we have \begin{equation}
\overline{A}_{k}(n)=\overline{B}_{k}(n).
\end{equation}
\end{thm}

\begin{thm}\label{lov2} For $k\geq 1$, let $\overline{D}_{k}(n)$ denote the number of overpartitions of $n$ of the
form $\lambda_1+\lambda_2+\cdots+\lambda_s$ such that $1$ cannot occur as a
non-overlined part, and  $\lambda_j-\lambda_{j+k-1}\geq1$ if $\lambda_j$ is
overlined and $\lambda_j-\lambda_{j+k-1}\geq2$ otherwise. Let $\overline{C}_{k}(n)$ denote
the number of overpartitions of $n$ whose non-overlined parts are
not congruent to $0,\pm1$ modulo $2k$. Then we have
\begin{equation}
\overline{C}_{k}(n)=\overline{D}_{k}(n).
\end{equation}
\end{thm}

 Chen, Sang and Shi \cite{chen} obtained  an overpartition analogue of the Rogers-Ramanujan-Gordon theorem in the general case.

\begin{thm}\label{thmlast1} For $k\geq i\geq 1$, let $P_{k,i}(n)$ denote the number of overpartitions of $n$ of the form $\lambda_1+\lambda_2+\cdots+\lambda_s$
such that part $1$  occurs as a non-overlined part at most $i-1$ times,
and  $\lambda_j-\lambda_{j+k-1}\geq1$ if $\lambda_j$ is overlined and
$\lambda_j-\lambda_{j+k-1}\geq2$ otherwise. For $k> i\geq 1$, let $Q_{k,i}(n)$ denote the number of
overpartitions of $n$ whose non-overlined parts are not congruent to
$0,\pm i$ modulo $2k$ and let $Q_{k,k}(n)$ denote the number of overpartitions of $n$ with parts not divisible by $k$.
 Then we have
 \begin{equation}
 P_{k,i}(n)=Q_{k,i}(n).
 \end{equation}
 \end{thm}

As an overpartition analogue of Bressoud's theorem for the case $i=1$,
Corteel,  Lovejoy, and  Mallet \cite{CLM} obtained the following relation.

\begin{thm}\label{CLM}For $k\geq 1$, let $\overline{A}_k^3(n)$ denote the number of overpartitions whose non-overlined parts are
not congruent to $0,\pm1$ modulo $2k-1$. Let $\overline{B}_k^3(n)$ denote the number of overpartitions $\lambda$ of $n$ such that
 (i) $f_1(\lambda)=0$, (ii) $f_l(\lambda)+f_{\overline{l}}(\lambda)+f_{l+1}(\lambda)\leq  k-1$, and (iii) if $f_l(\lambda)+f_{\overline{l}}(\lambda)+f_{l+1}(\lambda)=k-1$, then $lf_l(\lambda)+lf_{\overline{l}}(\lambda)+(l+1)f_{l+1}(\lambda)\equiv V_{\lambda}(l)\; (\rm{mod}\;2)$. Then we have \begin{equation}
 \overline{A}_k^3(n)=\overline{B}_k^3(n).
 \end{equation}
\end{thm}

In this paper, we shall give an overpartition analogue of the Bressoud's theorem
in the general case.

\section{The Main Result}

The main result of this paper can be stated as follows.

\begin{thm}\label{thm1}For $k\geq i\geq 1$, let $D_{k,i}(n)$ denote the number of overpartitions of $n$ of the form $\lambda=\lambda_1+\lambda_2+\cdots+\lambda_s$ such that
\begin{itemize}
  \item[(i)] $f_1(\lambda)\leq i-1$;
  \item[(ii)] $f_l(\lambda)+f_{\overline{l}}(\lambda)+f_{l+1}(\lambda)\leq k-1$; and
  \item[(iii)]  if $f_l(\lambda)+f_{\overline{l}}(\lambda)+f_{l+1}(\lambda)=k-1$, then $lf_l(\lambda)+lf_{\overline{l}}(\lambda)+(l+1)f_{l+1}(\lambda)\equiv V_{\lambda}(l)+i-1\; (\rm{mod} \; 2)$.
      \end{itemize}
      Let $C_{k,i}(n)$ denote the number of
overpartitions of $n$ whose nonoverlined parts are not congruent to
$0,\pm i$ modulo $2k-1$. Then we have
\begin{equation}
C_{k,i}(n)=D_{k,i}(n).
\end{equation}
 \end{thm}

In stead of using the function $\widetilde{J}_{k,i}(a; x; q)$ as in the proof of Theorem \ref{CLM} given  by Corteel,  Lovejoy, and  Mallet \cite{CLM}, we find that the
function $\widetilde{H}_{k,i}(a;x;q)$, also introduced by  Corteel,  Lovejoy, and  Mallet \cite{CLM}, is related to the generating functions of the numbers $C_{k,i}(n)$ and
$D_{k,i}(n)$.
Recall that
\begin{equation}\widetilde{J}_{k,i}(a;x;q)=\widetilde{H}_{k,i}(a;xq;q)+axq\widetilde{H}_{k,i-1}(a;xq;q),
\end{equation}
where
\begin{equation}\widetilde{H}_{k,i}(a;x;q)=\sum_{n\geq0}
\frac{(-a)^nq^{kn^2-{n\choose2}+n-in}x^{(k-1)n}(1-x^iq^{2ni})(-x,-1/a)_n(-axq^{n+1})_\infty}
{(q^2;q^2)_n(xq^n)_\infty}.
\end{equation}
It should be noticed that the function  $\widetilde{J}_{k,i}(a;x;q)$ can be expressed
as
$F_{1,k,i}(-q,\infty;-1/a;q)$  in the notation of  Bressoud \cite{bre80}, and the function $(-q)_\infty\widetilde{H}_{k,i}(a;x;q)$
can be written as  $H_{k,i}(-1/a,-x; x; q)_2$ in the notation of Andrews \cite{and68}.

Let $\widetilde{B}_k^3(m,n)$ denote the number of overpartitions enumerated by $\widetilde{B}_k^3(n)$ with exactly $m$ parts. Corteel,  Lovejoy and  Mallet  \cite{CLM}
have shown that
 the coefficients of $x^mq^n$ in $\widetilde{J}_{k,1}(1/q; x; q)$ and  $\widetilde{B}_k^3(m,n)$ satisfy the same recurrence relation with the same initial values.
Moreover, they   proved that the generating function of $\widetilde{B}_k^3(m,n)$ also equals $\widetilde{J}_{k,1}(1/q; x; q)$, that is,
\begin{equation}
\sum_{m,n\geq 0}\overline{B}_k^3(m,n)x^mq^n=\widetilde{J}_{k,1}(-1/q;x;q).
\end{equation}

Setting $a=-1/q$, $x=1$ and using the Jacobi's triple product identity,
 the function $\widetilde{J}_{k,i}(a;x;q)$ can be expressed as an infinite product, namely,
\[\widetilde{J}_{k,1}(-1/q;1;q)=
\frac{(q,q^{2k-2},q^{2k-1};q^{2k-1})_{\infty}(-q)_{\infty}}{(q)_{\infty}}.\]
Clearly, this is the generating function for $\overline{A}_k^3(n)$. Thus we have $\overline{A}_k^3(n)=\overline{B}_k^3(n)$.

However, the proof of Corteel,  Lovejoy and  Mallet does not seem to apply to the
general case, since $\widetilde{J}_{k,i}(-1/q;x;q)$ cannot be expressed as an infinite product for $i\geq 2$.
Our idea goes as follows.
For $C_{k,i}(n)$, we shall show that the generating function for $C_{k,i}(n)$ can be expressed in terms of $\widetilde{H}_{k,i}(a;x;q)$ with $a=-1/q$ and $x=q$. For $D_{k,i}(n)$, let $D_{k,i}(m,n)$ denote the number of overpartitions enumerated by $D_{k,i}(n)$ with exactly $m$ parts. We find a  combinatorial interpretation of $D_{k,i}(m,n)-D_{k,i-1}(m,n)$ from which we can derive a recurrence relation for $D_{k,i}(m,n)$.
Furthermore, we see that the recurrence relation and initial values of $D_{k,i}(m,n)$ coincide with the recurrence relation and  the initial values  of the coefficients of $x^mq^n$ in $\widetilde{H}_{k,i}(-1/q;xq;q)$. Thus we reach the conclusion that the generating function of $D_{k,i}(m,n)$ equals $\widetilde{H}_{k,i}(-1/q;xq;q)$. Setting $x=1$, we  deduce that the generating function of
$D_{k,i}(n)$ equals  the generating function of  $C_{k,i}(n)$.

For convenience, we write $W_{k,i}(x;q)$ for $\widetilde{H}_{k,i}(-1/q;xq;q)$, that is,
\begin{equation}W_{k,i}(x;q) =\sum_{n\geq 0}\frac{(-1)^nq^{(2k-1){n+1\choose2}-in}x^{(k-1)n}(1-x^iq^{(2n+1)i})(-xq)_\infty}{(q)_n(xq^{n+1})_\infty}.
\end{equation}

Recall that Andrews found the following recurrence relation for $H_{k,i}(a;x;q)$
\begin{equation}\label{hki}H_{k,i}(a;x;q)-H_{k,i-1}(a;x;q)
=x^{i-1}H_{k,k-i+1}(a;xq;q)-ax^iqH_{k,k-i}(a;xq;q).
\end{equation}
A recurrence relation for  $W_{k,i}(x;q)$ is given below.

\begin{thm}\label{WkiR}For $k\geq i\geq 1$, we have
\begin{equation}W_{k,i}(x;q)-W_{k,i-1}(x;q)=(1+xq)(xq)^{i-1}W_{k,k-i}(xq;q).
\end{equation}
\end{thm}

\pf  Since \[q^{-in}-x^iq^{(n+1)i}-q^{(-i+1)n}+x^{i-1}q^{(n+1)(i-1)}
=q^{-in}(1-q^n)+x^{i-1}q^{(n+1)(i-1)}(1-xq^{n+1}),\]
it can be checked that $W_{k,i}(x;q)-W_{k,i-1}(x;q)$ can be written as
\begin{equation}\label{ss}
\sum_{n=1}^{\infty}q^{-in}\frac{(-1)^nx^{(k-1)n}q^{(2k-1){n+1 \choose2}}(-xq)_{\infty}}
{(q)_{n-1}(xq^{n+1})_{\infty}}+\sum_{n=0}^{\infty}(xq^{n+1})^{i-1}
\frac{(-1)^nx^{(k-1)n}q^{(2k-1){n+1 \choose2}}(-xq)_{\infty}}
{(q)_n(xq^{n+2})_{\infty}}.
\end{equation}
Now, replacing $n$ with $n+1$, the first sum in (\ref{ss}) can be expressed as
\begin{equation} \label{ss2}
\sum_{n=0}^{\infty}q^{-i(n+1)}\frac{(-1)^{(n+1)}x^{(k-1)(n+1)}q^{(2k-1){n+2 \choose2}}(-xq)_{\infty}}
{(q)_{n}(xq^{n+2})_{\infty}}.
\end{equation}
Hence  $W_{k,i}(x;q)-W_{k,i-1}(x;q)$ equals
\begin{align*}&-(xq)^{i-1}\sum_{n=0}^{\infty}\frac{(-1)^n(xq)^{(k-1)n}q^{(2k-1){n+1 \choose2}}x^{k-i}q^{(2k-1)(n+1)-in-2i+1-(k-1)n}(-xq)_{\infty}}
{(q)_{n}(xq^{n+2})_{\infty}}\\&\quad +(xq)^{i-1}\sum_{n=0}^{\infty}
\frac{(-1)^n(xq)^{(k-1)n}q^{(2k-1){n+1 \choose2}+(i-1)n-(k-1)n}(-xq)_{\infty}}{(q)_n(xq^{n+2})_{\infty}}
\\& \qquad = (1+xq)(xq)^{i-1}\sum_{n\geq0}\frac{(-1)^n(xq)^{(k-1)n}q^{(2k-1){n+1 \choose2}-(k-i)n}(1-x^{k-i}q^{(2n+2)(k-i)})(-xq^2)_\infty}{(q)_n(xq^{n+2})_{\infty}}
\\& \qquad = (1+xq)(xq)^{i-1}W_{k,k-i}(xq),
\end{align*}
as desired. \qed

The following relation can be considered as a combinatorial interpretation
of  $D_{k,i}(m,n)-D_{k,i-1}(m,n)$.

\begin{thm}For $k\geq i\geq 1$ and for $m,n\geq0$, let $S_{k,i}(m,n)$ denote the set of the overpartitions enumerated by $D_{k,i}(m,n)$ with exactly one overlined part $1$ and exactly $i-1$ nonoverlined part $1$. Let $T_{k,i}(m,n)$ denote  the set of the overpartitions enumerated by $D_{k,i}(m,n)$ with exactly one overlined part $1$ and exactly $i-1$ nonoverlined part $1$. Let $Q_{k,i}(m,n)$ denote the number of overpartitions in $S_{k,i}(m,n)$ and let $R_{k,i}(m,n)$ denote the number of overpartitions in $T_{k,i}(m,n)$. Then we have

\begin{equation}\label{DQR}
 D_{k,i}(m,n)-D_{k,i-1}(m,n)=Q_{k,i}(m,n)+R_{k,i}(m,n).
\end{equation}
\end{thm}
\pf Let $U_{k,i}(m,n)$ denote the set of overpartitions enumerated by $D_{k,i}(n)$ with exactly $m$ parts. By the definition of $D_{k,i}(m,n)$ and $D_{k,i-1}(m,n)$, it can easily seen that  $U_{k,i-1}(m,n)$ is not contained in $U_{k,i}(m,n)$.
  To compute $D_{k,i}(m,n)-D_{k,i-1}(m,n)$,  we wish to construct
     an injection $\varphi$ from   overpartitions in $U_{k,i-1}(m,n)$ to overpartitions $U_{k,i}(m,n)$.  From the characterization of the images of this map, we obtain the
     relation (\ref{DQR}).

     Let $\lambda$ be an overpartition in $U_{k,i-1}(m,n)$.
  If there exists an overlined part  of  $\lambda$ with the smallest underlying part,  then we switch this overlined   part to a nonoverlined part, otherwise we choose a smallest nonoverlined  part and switch it to an overlined part. Let $\lambda'$ denote the resulting overpartition.
   It can be checked that this map is an injection. It is not difficult to verify
     that $\lambda'\in U_{k,i}(m,n)$.
   Hence  the number
      $D_{k,i}(m,n)-D_{k,i-1}(m,n)$ can be interpreted as
      the number of  overpartitions in
      $U_{k,i}(m,n)$ which cannot be obtained by using the above map.

      By the construction of the map $\varphi$,  we may generate
       all the overpartitions in $U_{k,i}(m,n)$ with no overlined part equal to $1$ and all the overpartitions in $U_{k,i}(m,n)$ with an overlined $1$ and with at most $i-3$ nonoverlined part $1$. Therefore, $D_{k,i}(m,n)-D_{k,i-1}(m,n)$ is exactly  the number of overpartitions in $U_{k,i}(m,n)$ with exactly one overlined part $1$ such that
       the  nonoverlined part $1$ appears either $i-1$ or $i-2$ times. This completes the proof. \qed

\begin{thm}For $k\geq i\geq1$,and $m,n \geq 0$, we have
\begin{equation} \label{Qki}
 Q_{k,i}(m,n)=D_{k,k-i}(m-i,n-m).
 \end{equation}
\end{thm}

\pf
We shall  define a bijection $\phi$ from $S_{k,i}(m,n)$ to
$U_{k,k-i}(m-i,n-m)$ which implies
(\ref{Qki}).
 Let  $\lambda$ be  an overpartition in $S_{k,i}(m,n)$, the map $\phi$ is defined
  as follows.

 \noindent Step 1. Remove all the $i-1$ parts with underlying part  $1$.

 \noindent Step 2. Subtract $1$ from each part.

Clearly, the resulting overpartition $\lambda'$ is an overpartition of $n-m$ with $m-i$ parts.  Moreover, we
 claim that $\lambda'\in U_{k,k-i}(m-i,n-m)$.

 We first show that $f_{1}(\lambda')\leq k-i-1$. By the construction of $\phi$,
 it is easy to see that $f_{1}(\lambda')=f_2(\lambda)$ and $f_1(\lambda)=i-1$. From  the condition (ii) in the theorem, that is,  $f_1(\lambda)+f_{\overline{1}}(\lambda)+f_2(\lambda)\leq k-1$, we find that $f_2(\lambda)\leq k-1-i$.

We still need to verify that if there is an integer $l$ such that \begin{equation}\label{ll}f_{l}(\lambda')+f_{\overline{l}}(\lambda')+f_{l+1}(\lambda')=k-1,
\end{equation}
then we have
\begin{equation}\label{lf}
lf_{l}(\lambda')+lf_{\overline{l}}(\lambda')+(l+1)f_{l+1}(\lambda')\equiv V_{\lambda'}(l)+k-i-1\;(\rm{mod}\;2).\end{equation}
By the construction of $\phi$, it is easily checked that  (\ref{ll}) implies
  \[ f_{l+1}(\lambda)+f_{\overline{l+1}}(\lambda)+f_{l+2}(\lambda)=k-1.
\]
Since $\lambda \in S_{k,i}(m,n)$, we have \[(l+1)f_{l+1}(\lambda)+(l+1)f_{\overline{l+1}}(\lambda)+(l+2)f_{l+2}(\lambda)\equiv V_{\lambda}(l+1)+i-1\;(\rm{mod}\;2).\]
Clearly, $f_{t}(\lambda')=f_{t+1}(\lambda)$ and $f_{\overline{t}}(\lambda')=f_{\overline{t+1}}(\lambda)$ for any $t\geq 1$.
Thus we deduce that
\[lf_{l}(\lambda')+lf_{\overline{l}}(\lambda')+(l+1)f_{l-1}(\lambda')\equiv V_{\lambda}(l+1)+i-1-(k-1)\;(\rm{mod}\;2).\]
Again, by the construction of $\phi$, we find $V_{\lambda}(l+1)=V_{\lambda'}(l)+1$.
So we arrive at  relation (\ref{lf}), which implies   $\lambda'\in U_{k,k-i}(m-i,n-m)$.

It is not difficult to verify that the above construction is reversible, that is,
from any overpartition in $U_{k,k-i}(m-i,n-m)$, we can recover an overpartition in $S_{k,i}(m,n)$. This completes the
proof.   \qed

\begin{thm}For $k\geq i\geq1$,and $m,n \geq 0$, we have
\begin{equation}\label{RD}
 R_{k,i}(m,n)=D_{k,k-i}(m-i+1,n-m).
 \end{equation}
\end{thm}

\pf We proceed to give a bijection $\chi$ from $T_{k,i}(m,n)$ to $U_{k,k-i}(m-i+1,n-m)$.
Let $\lambda$ be an overpartition in $T_{k,i}(m,n)$, the map $\chi$ is defined as follows.

\noindent Step 1. Remove all $i-1$ parts equal to $1$.

\noindent Step 2. Subtract $1$ from each part.

Clearly, the resulting overpartition $\lambda'$ is an overpartition of $n-m$ with $m-i+1$ parts.  We shall show that $\lambda'\in U_{k,k-i}(m-i+1,n-m)$.

We first verify that $f_{1}(\lambda')\leq k-i-1$. It is obvious that $f_{1}(\lambda')=f_{2}(\lambda)$. So it suffices to prove that $f_{2}(\lambda)\leq k-i-1$.
Since $\lambda\in T_{k,i}(m,n)$, we have $f_{1}(\lambda)=i-2$, $f_{\overline{1}}(\lambda)=1$ and \begin{equation}\label{f1}
 f_1(\lambda)+f_{\overline{1}}(\lambda)+f_2(\lambda)\leq k-1.\end{equation}
It follows that  $f_2(\lambda)\leq k-i$.

It remains to show that the nonoverlined part $2$ cannot occur $k-i$ times. Assume that  $f_2(\lambda)= k-i$. Then the equality in (\ref{f1}) holds, that is,
 \[f_1(\lambda)+f_{\overline{1}}(\lambda)+f_2(\lambda)=k-1.\]
We wish to derive a contradiction to the condition (iii) in Theorem \ref{thm1}. By the facts  $f_1(\lambda)=i-2$, $f_{\overline{1}}(\lambda)=1$,  we find
\begin{equation} \label{1f1}
1f_{1}(\lambda)+1f_{\overline{1}}(\lambda)+2f_{2}(\lambda)=2k-i-1.\end{equation}
Since $V_{\lambda}(1)=1$, from (\ref{1f1}) it follows that
 \[1f_{1}(\lambda)+1f_{\overline{1}}(\lambda)+2f_{2}(\lambda)\not \equiv V_{\lambda}(1)+i-1\;(\rm{mod}\;2),\] a contradiction. Thus we reach the
 conclusion that the nonoverlined part $2$ occurs at most $k-i-1$ times in $\lambda$, or
  equivalently, the nonoverlined part $1$  occurs at most $k-i-1$ times in $\lambda'$.

  Next, we check condition (ii) for $\lambda'$. For any $l\geq 1$, we see that \begin{equation}\label{l}f_{l+1}(\lambda)=f_{l}(\lambda')\quad \mbox{and} \quad f_{\overline{l+1}}(\lambda)=f_{\overline{l}}(\lambda').\end{equation}  From condition (ii)
  for $\lambda$, we get
\[f_{l}(\lambda')+f_{\overline{l}}(\lambda')+f_{l+1}(\lambda')\leq k-1.\]

Finally, we proceed to verify that if there is an  integer $l$ such that \begin{equation}\label{fl2}f_{l}(\lambda')
+f_{\overline{l}}(\lambda')+f_{l+1}(\lambda')=k-1,\end{equation}
then we have
\begin{equation}\label{ll2}lf_{l}(\lambda')+lf_{\overline{l}}(\lambda')+(l+1)f_{l+1}(\lambda')\equiv V_{\lambda'}(l)+k-i-1\;(\rm{mod}\;2).\end{equation}

Notice that (\ref{fl2}) implies \begin{equation}\label{ll4}f_{l+1}(\lambda)+f_{\overline{l+1}}(\lambda)+f_{l+2}(\lambda)=k-1.\end{equation}
Since $\lambda \in T_{k,i}(m,n)$, by  condition (iii) for $\lambda$, we have  \begin{equation}\label{l5}(l+1)f_{l+1}
(\lambda)+(l+1)f_{\overline{l+1}}(\lambda)+(l+2)f_{l+2}(\lambda)\equiv V_{\lambda}(l+1)+i-1\;(\rm{mod}\;2).\end{equation}
 Substituting \eqref{l} into \eqref{l5}, we obtain  \begin{equation}\label{l6} lf_{l}(\lambda')+lf_{\overline{l}}(\lambda')+(l+1)f_{l+1}(\lambda')\equiv V_{\lambda}(l+1)+i-1-(k-1)\;(\rm{mod}\;2).\end{equation}
 Observing that $V_{\lambda}(l+1)=V_{\lambda'}(l)+1$, (\ref{l6}) can be rewritten as (\ref{ll2}). This leads to the conclusion that  $\lambda'\in U_{k,k-i}(m-i+1,n-m)$.

  It is routine to verify that the above procedure is reversible, that is,
  from  any overpartition in $U_{k,k-i}(m-i+1,n-m)$, one can recover an overpartition
   in $T_{k,i}(m,n)$. This completes the proof. \qed

By relations \eqref{DQR}, \eqref{Qki} and \eqref{RD},  we obtain a recurrence relation of $D_{k,i}(m,n)$.

\begin{thm}\label{Dki}For $k\geq i\geq 1$ and for $m,n\geq 0$, we have
\begin{equation}\label{dki}D_{k,i}(m,n)-D_{k,i-1}(m,n)=D_{k,k-i}(m-i,n-m)+D_{k,k-i}(m-i+1,n-m).
\end{equation}
\end{thm}

By Theorem \ref{WkiR} and Theorem \ref{Dki}, we obtain a combinatorial interpretation of $W_{k,i}(x;q)$ in terms of overpartitions.

\begin{thm}\label{DH}For $k\geq i\geq 1$, we have
\begin{equation}\label{DW}W_{k,i}(x;q)=\sum_{m,n\geq 0}D_{k,i}(m,n)x^mq^n.
\end{equation}
\end{thm}

\pf For $m, n\geq 0$ and for $k\geq i \geq 1$, let $W_{k,i}(m,n)$ denote
the coefficient of $x^mq^n$ in $W_{k,i}(x;q)$, that is,
\begin{equation}\label{wki}W_{k,i}(x;q)=\sum_{m,n\geq 0}^{\infty}W_{k,i}(m,n)x^mq^n.
\end{equation}
We proceed to show that $D_{k,i}(m,n)$ and $W_{k,i}(m,n)$ satisfy the
same recurrence relations with the same initial values.

 Clearly, we have
 $W_{k,i}(0,0)=1$ for $ k\geq i\geq 1$ and
 $W_{k,0}(m,n)=0$ for $k\geq 1, m,n\geq 0$. Moreover, we assume that
 $W_{k,i}(m,n)=0$ if $m$ or $n$ is zero but not both.
By Theorem \ref{WkiR}, we find that
\begin{equation}
\label{j4}W_{k,i}(m,n)-W_{k,i-1}(m,n)=W_{k,k-i}(m-i,n-m)+W_{k,k-i}(m-i+1,n-m),
\end{equation}
which is the same recurrence relation as $D_{k,i}(m,n)$ as given in Theorem \ref{Dki}.

  It is clear that $D_{k,i}(0,0)=1$ for $ k\geq i\geq 1$ and
 $D_{k,0}(m,n)=0$ for $k\geq 1, m,n\geq 0$. Moreover,
 $D_{k,i}(m,n)=0$ if $m$ or $n$ is zero but not both. Now, we see that $D_{k,i}(m,n)$
 and $W_{k,i}(m,n)$ have the same recurrence relation and the same initial values.
 This completes the proof.  \qed

We are now ready to finish the proof of Theorem \ref{thm1}.

\noindent {\it Proof of Theorem \ref{thm1}.} Setting $x=1$ in \eqref{DW},  we find that the generating function for $D_{k,i}(n)$ equals $W_{k,i}(1;q)$. In other words,
\begin{equation}\label{DD}\sum_{n\geq0}D_{k,i}(n)q^n=\sum_{n=0}^{\infty}\frac{(-1)^nq^{(2k-1){n+1 \choose2}-in}(1-q^{(2n+1)i})(-q)_{\infty}}
{(q)_n(q^{n+1})_{\infty}}.
\end{equation}
 The right hand side of  \eqref{DD} can be expressed as \begin{equation}\label{D3}\frac{(-q)_{\infty}}{(q)_{\infty}}\sum_{n=0}^{\infty}(-1)^nq^{(2k-1){n+1 \choose2}-in}+\frac{(-q)_{\infty}}{(q)_{\infty}}\sum_{n=0}^{\infty}(-1)^nq^{(2k-1){n+1 \choose2}+i(n+1)}.\end{equation}
 By substituting $n$ with $-(n+1)$ in the second sum of (\ref{D3}), we get \[\frac{(-q)_{\infty}}{(q)_{\infty}}
\sum_{n=-\infty}^{\infty}(-1)^nq^{(2k-1){n+1 \choose2}-in}.\]
In view of Jacobi's triple product identity,  we obtain
\begin{equation}\label{Dgen}\sum_{n\geq0}D_{k,i}(n)q^n=
\frac{(q^i,q^{2k-1-i},q^{2k-1};q^{2k-1})_{\infty}(-q)_{\infty}}{(q)_{\infty}}.\end{equation}

By the definition of $C_{k,i}(n)$, it is easily seen that
\begin{equation}\label{Cgen}\sum_{n=0}^{\infty}C_{k,i}(n)q^n=
\frac{(q^i,q^{2k-1-i},q^{2k-1};q^{2k-1})_{\infty}(-q)_{\infty}}{(q)_{\infty}}.
\end{equation}
Comparing \eqref{Dgen} and \eqref{Cgen} we deduce that
$C_{k,i}(n)=D_{k,i}(n)$ for  $k\geq i\geq 1$ and $n\geq 0$.
This completes the proof. \qed

\vspace{0.5cm}
 \noindent{\bf Acknowledgments.}  This work was supported by  the 973
Project, the PCSIRT Project of the Ministry of Education,  and the
National Science Foundation of China.

\end{document}